\documentclass{elsart}
\usepackage{amsfonts}
\newcommand{\eps}{\varepsilon}
\newcommand{\C}{\mathbb{C}}
\newcommand{\R}{\mathbb{R}}
\begin{document}
\begin{frontmatter}
\title{%
Maximal operators on variable Lebesgue\\
spaces with weights related to oscillations\\
of Carleson curves}
\author{Alexei Yu. Karlovich}
\thanks{This work is partially supported by the grant of
F.C.T. (Portugal) FCT/ FEDER/POCTI/MAT/59972/2004}
\address{%%
Departamento de Matem\'atica,
Faculdade de Ci\^encias e Tecnologia,
Universidade Nova de Lisboa,
Quinta da Torre,
2829--516 Caparica,
Portugal}
\ead{oyk@fct.unl.pt}

\begin{abstract}
We prove sufficient conditions for the boundedness of the maximal operator on
variable Lebesgue spaces with weights $\varphi_{t,\gamma}(\tau)=|(\tau-t)^\gamma|$,
where $\gamma$ is a complex number, over arbitrary Carleson curves. If the curve has
different spirality indices at the point $t$ and $\gamma$ is not real, then
$\varphi_{t,\gamma}$ is an oscillating weight lying beyond the class of radial
oscillating weights considered recently by V. Kokilashvili, N. Samko, and S. Samko.
\end{abstract}
\end{frontmatter}
%%%%%%%%%%%%%%%%%%%%%%%%%%%%%%%%%%%%%%%%%%%%%%%%%%%%%%%%%%%%%%%%%%%%%%%%%%%%%%%%%%%%%%%%%%
\section{Introduction and main result}
Let $\Gamma$ be a rectifiable curve in the complex plane equipped with arc-length measure
$|d\tau|$. We suppose that $\Gamma$ is simple, that is, homeomorphic to a segment or to a
circle. A measurable function $w:\Gamma\to[0,\infty]$ is said to be a weight if it is
positive and finite almost everywhere. Let $p:\Gamma\to (1,\infty)$ be a continuous
function. A weighted variable Lebesgue space $L^{p(\cdot)}(\Gamma,w)$ is the set of all
measurable complex-valued functions $f$ on $\Gamma$ such that
\[
\int_\Gamma |f(\tau)w(\tau)/\lambda|^{p(\tau)}\,|d\tau|<\infty
\]
for some $\lambda=\lambda(f)>0$. It is a Banach space when equipped with the Luxemburg-Nakano
norm
\[
\|f\|_{L^{p(\cdot)}(\Gamma,w)}=\inf\left\{
\lambda>0\ :\ \int_\Gamma |f(\tau)w(\tau)/\lambda|^{p(\tau)}\,|d\tau|\le 1
\right\}.
\]
It is clear that $L^{p(\cdot)}(\Gamma,w)$ coincides with the standard Lebesgue space
whenever $p$ is constant. It is a partial case of so-called Musielak-Orlicz spaces
(see \cite{Musielak83,MO59}).

Two weights $w_1$ and $w_2$ on $\Gamma$ are said to be equivalent
if there is a bounded and bounded away from zero function $f$ on $\Gamma$
such that $w_1=fw_2$. It is easy to see that $L^{p(\cdot)}(\Gamma,w_1)$ and
$L^{p(\cdot)}(\Gamma,w_2)$ are isomorphic whenever $w_1$ and $w_2$ are equivalent.

A curve $\Gamma$ is said to be Carleson (or Ahlfors-David regular) if
\[
C_\Gamma:=\sup_{t\in\Gamma}\sup_{\eps>0}\frac{|\Gamma(t,\eps)|}{\eps}<\infty
\]
where $\Gamma(t,\eps):=\{\tau\in\Gamma:|\tau-t|<\eps\}$ is the portion of the curve
in the disk centered at $t$ of radius $\eps$ and $|\Omega|$ denotes the measure of
a measurable set $\Omega\subset\Gamma$. We are interested in the boundedness conditions
for the maximal operator
\[
(Mf)(t):=\sup_{\eps>0}\frac{1}{|\Gamma(t,\eps)|}\int_{\Gamma(t,\eps)}|f(\tau)|\,|d\tau|
\quad
(t\in\Gamma)
\]
on weighted variable Lebesgue spaces. This operator is one of the main players in
harmonic analysis. It is closely related to the Cauchy singular integral operator
\[
(Sf)(t):=\lim_{\eps\to 0}\frac{1}{\pi i}\int_{\Gamma\setminus\Gamma(t,\eps)}
\frac{f(\tau)}{\tau-t}\,d\tau
\quad(t\in\Gamma).
\]
The boundedness of both operators on standard weighted Lebesgue spaces is well understood
(see e.g. \cite{BK97,GGKK98,Dynkin91,Stein93}).
If $T$ is one of the operators $M$ or $S$ and $1<p<\infty$, then $T$ is bounded on
$L^p(\Gamma,w)$ if and only if $w$ is a Muckenhoupt weight, $w\in A_p(\Gamma)$, that is,
\[
\sup_{t\in\Gamma}\sup_{\eps>0}
\left(\frac{1}{\eps}\int_{\Gamma(t,\eps)}w^p(\tau)\,|d\tau|\right)^{1/p}
\left(\frac{1}{\eps}\int_{\Gamma(t,\eps)}w^{-q}(\tau)\,|d\tau|\right)^{1/q}<\infty
\]
where $1/p+1/q=1$. By H\"older's inequality, if $w$ is a Muckenhoupt weight, then
$\Gamma$ is a Carleson curve.

Let us define the weight we are interested in.
Fix $t\in\Gamma$ and consider the function $\eta_{t}:\Gamma\setminus\{t\}\to(0,\infty)$
defined by
\[
\eta_{t}(\tau):=e^{-\arg(\tau-t)},
\]
where $\arg(\tau-t)$ denotes any continuous branch of the argument on $\Gamma\setminus\{t\}$.
For every $\gamma\in\C$, put
%%%%
\[
\varphi_{t,\gamma}(\tau):=|(\tau-t)^\gamma|=
|\tau-t|^{{\rm Re}\,\gamma}\eta_t(\tau)^{{\rm Im}\,\gamma}
\quad (\tau\in\Gamma\setminus\{t\}).
\]
In the Fredholm theory of singular integral operators with piecewise continuous
coefficients on $L^{p(\cdot)}(\Gamma)$ (without weights!) it is important to know
for which values of $\gamma$ the operator $S$ is bounded on
$L^{p(\cdot)}(\Gamma,\varphi_{t,\gamma})$ (see e.g.
\cite{Karlovich03-JIEA,Karlovich08-IWOTA} and also \cite{BK97}).
In fact, an attempt to answer this question is the our main motivation for this work.

The above question was completely studied for the case of standard Lebesgue
spaces by A.~B\"ottcher and Yu.~Karlovich \cite[Section~3.1]{BK97}.
To formulate their result explicitly, we need some definitions.
A function $\varrho:(0,\infty)\to(0,\infty]$ is called regular if it is bounded from above
in some open neighborhood of the point $1$. A function $\varrho:(0,\infty)\to(0,\infty]$
is said to be submultiplicative if $\varrho(x_1x_2)\le\varrho(x_1)\varrho(x_2)$
for all $x_1,x_2\in(0,\infty)$. A regular submultiplicative function is finite everywhere
and one can define
\[
\alpha(\varrho):=\sup_{x\in(0,1)}\frac{\log\varrho(x)}{\log x},
\quad
\beta(\varrho):=\inf_{x\in(1,\infty)}\frac{\log\varrho(x)}{\log x}.
\]
In this case, by \cite[Theorem~1.13]{BK97}, $-\infty<\alpha(\varrho)\le\beta(\varrho)<\infty$.
The numbers $\alpha(\varrho)$ and $\beta(\varrho)$ are called lower and upper indices of
$\varrho$, respectively.

For $t\in\Gamma$, put $d_{t}:=\max\limits_{\tau\in\Gamma}|\tau-t|$.
Following \cite[Section~1.5]{BK97}, for a continuous function $\psi:\Gamma\setminus\{t\}\to(0,\infty)$,
we define
\[
(W_{t}\psi)(x):=\left\{
\begin{array}{lll}
\displaystyle
\sup_{0<R\le d_{t}}\left(
\max_{|\tau-t|=xR}\psi(\tau)/\min_{|\tau-t|=R}\psi(\tau)
\right)
& \mbox{for} & x\in(0,1],\\
\displaystyle
\sup_{0<R\le d_{t}}\left(
\max_{|\tau-t|=R}\psi(\tau)/\min_{|\tau-t|=x^{-1}R}\psi(\tau)
\right)
& \mbox{for} & x\in[1,\infty).
\end{array}
\right.
\]
This function is submultiplicative in view of \cite[Lemma~1.15]{BK97}.
From \cite[Theorem~1.18]{BK97} it follows that $W_{t}\eta_{t}$ is regular for every
$t\in\Gamma$ whenever $\Gamma$ is a Carleson curve. Hence the lower and upper spirality
indices $\delta_{t}^-$ and $\delta_{t}^+$ at $t\in\Gamma$ are correctly defined by
\[
\delta_{t}^-:=\alpha(W_{t}\eta_{t}),
\quad
\delta_{t}^+:=\beta(W_{t}\eta_{t}).
\]
%%%%%%%%%%%%%%%%%%%%%%%%%%%%%%%%%%%%%%%%%%%%%%%%%%%%%%%%%%%%%%%%%%%%%%%%%%%%%%%%%%%
\begin{prop}\label{pr:Carleson}
\begin{enumerate}
\item[{\rm (a)}]
If $\Gamma$ is a piecewise smooth curve, then $\arg(\tau-t)=O(1)$
and $\delta_t^-=\delta_t^+=0$ for all $t\in\Gamma$.

\item[{\rm (b)}] If $\Gamma$ is a Carleson curve satisfying
%%%%
\begin{equation}\label{eq:logarithmic-Carleson}
\arg(\tau-t)=-\delta\log|\tau-t|+O(1)
\quad\mbox{as}\quad\tau\to t
\end{equation}
%%%%
at some $t\in\Gamma$ with some $\delta\in\R$, then $\delta_t^-=\delta_t^+=\delta$.

\item[{\rm (c)}] {\rm\textbf{(R. Seifullayev)}}
If $\Gamma$ is a Carleson curve, then
%%%
\begin{equation}\label{eq:general-Carleson}
\arg(\tau-t)=O(-\log|\tau-t|)
\quad\mbox{as}\quad\tau\to t
\end{equation}
%%%
for every $t\in\Gamma$.

\item[{\rm (d)}] {\rm\textbf{(A. B\"ottcher, Yu. Karlovich)}}
For any given real numbers $\alpha,\beta$ such that
\[
-\infty<\alpha<\beta<+\infty,
\]
there exists a Carleson curve $\Gamma$ such that $\delta_t^-=\alpha$ and $\delta_t^+=\beta$
at some point $t\in\Gamma$.
\end{enumerate}
\end{prop}
%%%%%%%%%%%%%%%%%%%%%%%%%%%%%%%%%%%%%%%%%%%%%%%%%%%%%%%%%%%%%%%%%%%%%%%%%%%%%%
Parts (a) and (b) are trivial, a proof of part (c) is in \cite[Theorem~1.10]{BK97},
and part (d) is proved in \cite[Proposition~1.21]{BK97}.
From \cite[Propistion~3.1]{BK97} it follows that $W_t\varphi_{t,\gamma}$ is regular for
every $\gamma\in\C$ and
\begin{eqnarray*}
\alpha(W_t\varphi_{t,\gamma})&=&{\rm Re}\,\gamma+
\min\{\delta_t^-{\rm Im}\,\gamma,\delta_t^+{\rm Im}\,\gamma\},
\\
\beta(W_t\varphi_{t,\gamma})&=&{\rm Re}\,\gamma+
\max\{\delta_t^-{\rm Im}\,\gamma,\delta_t^+{\rm Im}\,\gamma\}.
\end{eqnarray*}
These equalities in conjunction with \cite[Theorem~2.33]{BK97} yield the following.
%%%%%%%%%%%%%%%%%%%%%%%%%%%%%%%%%%%%%%%%%%%%%%%%%%%%%%%%%%%%%%%%%%%%%%%%%%%%%%
\begin{thm}[A.~B\"ottcher, Yu.~Karlovich]
\label{th:BK}
Let $\Gamma$ be a Carleson curve and $p\in(1,\infty)$ be constant. Suppose
$t\in\Gamma$ and $\gamma\in\C$. Then $\varphi_{t,\gamma}\in A_p(\Gamma)$ if and only if
\begin{eqnarray*}
0&<&\frac{1}{p}+{\rm Re}\,\gamma+\min\{\delta_t^-{\rm Im}\,\gamma,\delta_t^+{\rm Im}\,\gamma\}
\\
&\le&\frac{1}{p}+{\rm Re}\,\gamma+\min\{\delta_t^-{\rm Im}\,\gamma,\delta_t^+{\rm Im}\,\gamma\}<1.
\end{eqnarray*}
\end{thm}
%%%%%%%%%%%%%%%%%%%%%%%%%%%%%%%%%%%%%%%%%%%%%%%%%%%%%%%%%%%%%%%%%%%%%%%%%%%%%%

In the last decade many results from classical harmonic analysis for standard (weighted)
Lebesgue spaces were extended to the setting of (weighted) variable Lebesgue spaces
(see e.g. \cite{CFMP06,Diening04,KS04,KS08,Lerner05} and the references therein).
We recall only the most relevant result. Following \cite{Diening04,KS04}, we will
always suppose that $p:\Gamma\to(1,\infty)$ is a continuous function satisfying
the Dini-Lipschitz condition on $\Gamma$, that is, there exists a constant $C_p>0$
such that
%%%
\begin{equation}\label{eq:Dini-Lipschitz}
|p(\tau)-p(t)|\le\frac{C_p}{-\log|\tau-t|}
\end{equation}
%%%
for all $\tau,t\in\Gamma$ such that $|\tau-t|\le1/2$. For power weights one has the next
criterion. For simplicity, we formulate it in the case of one singularity only.
However, it is valid for power weights with a finite number of singularities
(see \cite{KPS06-OTAA,KS05-Simonenko}).
%%%%%%%%%%%%%%%%%%%%%%%%%%%%%%%%%%%%%%%%%%%%%%%%%%%%%%%%%%%%%%%%%%%%%%%%%%%%%%%%%%%%%%%%
\begin{thm}[V.~Kokilashvili, V.~Paatashvili, S.~Samko]
\label{th:KPS}
Let $\Gamma$ be a Carleson curve and $p:\Gamma\to(1,\infty)$ be a continuous function
satisfying the Dini-Lipschitz condition. For $t\in\Gamma$ and $\lambda\in\R$,
define the power weight $w(\tau):=|\tau-t|^\lambda$.
If $T$ is one of the operators $M$ or $S$, then $T$ is bounded on $L^{p(\cdot)}(\Gamma,w)$
if and only if
\[
0<\frac{1}{p(t)}+\lambda<1.
\]
\end{thm}
%%%%%%%%%%%%%%%%%%%%%%%%%%%%%%%%%%%%%%%%%%%%%%%%%%%%%%%%%%%%%%%%%%%%%%%%%%%%%%%%%%%%%%%%
Clearly, $\varphi_{t,\gamma}$ is equivalent to a power weight
$w(\tau)=|\tau-t|^\lambda$ if and only if $\gamma$ is real or $\Gamma$
satisfies (\ref{eq:logarithmic-Carleson}) at $t$. Hence, Theorem~\ref{th:KPS}
is not applicable to the weight $\eta_t$ for Carleson curves with $\delta_t^-<\delta_t^+$.

The sufficiency portion of Theorem~\ref{th:KPS} has been extended recently to
the case of radial oscillating weights (see \cite{KSS07-JFSA} for $M$
and \cite{KSS07-MN} for $S$).
In the case of one singularity these weights have the form $w(\tau)=f(|\tau-t|)$
where $t\in\Gamma$ is fixed and $f:(0,{\rm diam}(\Gamma)]\to(0,\infty)$
is some continuous function with additional regularity properties. It is clear that
the weight $\eta_t$ is not of this form. Thus, in general, weights considered
in this paper lie beyond the class of radial oscillating weights.
%%%%%%%%%%%%%%%%%%%%%%%%%%%%%%%%%%%%%%%%%%%%%%%%%%%%%%%%%%%%%%%%%%%%%%%%%%%%%%%%%%%%%%%%
\begin{thm}[Main result]
\label{th:main}
Let $\Gamma$ be a Carleson curve and $p:\Gamma\to(1,\infty)$ be a continuous function
satisfying the Dini-Lipschitz condition. If $t\in\Gamma$, $\gamma\in\C$, and
\begin{eqnarray}
0&<&
\frac{1}{p(t)}+{\rm Re}\,\gamma+
\min\{\delta_t^-{\rm Im}\,\gamma,\delta_t^+{\rm Im}\,\gamma\}
\nonumber
\\
&\le&
\frac{1}{p(t)}+{\rm Re}\,\gamma+
\max\{\delta_t^-{\rm Im}\,\gamma,\delta_t^+{\rm Im}\,\gamma\}<1,
\label{eq:main-conditions}
\end{eqnarray}
then $M$ is bounded on $L^{p(\cdot)}(\Gamma,\varphi_{t,\gamma})$.
\end{thm}
%%%%%%%%%%%%%%%%%%%%%%%%%%%%%%%%%%%%%%%%%%%%%%%%%%%%%%%%%%%%%%%%%%%%%%%%%%%%%%%%%%%%%%%%
We conjecture that Theorem~\ref{th:main} is true with $M$ replaced by $S$ and that a check
of the proof of \cite[Theorem~4.3]{KSS07-MN} will indicate the modifications needed to
obtain the desired result. We also conjecture that inequalities (\ref{eq:main-conditions})
are necessary for the boundedness of $M$ and $S$ on $L^{p(\cdot)}(\Gamma,\varphi_{t,\gamma})$.
To support the second conjecture, note that arguing as in \cite{Karlovich08-JFSA}, one can
show that if $S$ is bounded on $L^{p(\cdot)}(\Gamma,\varphi_{t,\gamma})$, then
\begin{eqnarray*}
0&\le&
\frac{1}{p(t)}+{\rm Re}\,\gamma+
\min\{\delta_t^-{\rm Im}\,\gamma,\delta_t^+{\rm Im}\,\gamma\}
\nonumber
\\
&\le&
\frac{1}{p(t)}+{\rm Re}\,\gamma+
\max\{\delta_t^-{\rm Im}\,\gamma,\delta_t^+{\rm Im}\,\gamma\}\le 1.
\end{eqnarray*}

The paper is organized as follows. In Section~\ref{sec:Muckenhoupt-ersatz}
we formulate a sufficient condition for the boundedness of $M$ on $L^{p(\cdot)}(\Gamma,w)$
involving the classical Muckenhoupt condition. Further we apply it to the case of the weight
$\varphi_{t,\gamma}$. In Section~\ref{sec:estimate} we estimate a weight $w$ with the
only singularity at $t$ by power weights with exponents $\alpha(W_tw)-\eps$
and $\beta(W_tw)+\eps$ where $\eps$ is small enough. Section~\ref{sec:proof} contains the
proof of Theorem~\ref{th:main}. Here we follow an idea from \cite{KSS07-JFSA} and represent
the weighted maximal operator as the sum of four maximal operators. The first operator is the maximal
operator over a small arc containing the singularity of the weight $\varphi_{t,\gamma}$.
Its boundedness follows from the results of Section~\ref{sec:Muckenhoupt-ersatz}.
The second and third maximal operators are estimated by maximal operators with power weights
with exponents $\alpha(W_t\varphi_{t,\gamma})-\eps$ and $\beta(W_t\varphi_{t,\gamma})+\eps$
by using the results of Section~\ref{sec:estimate}. The boundedness of the latter operators
follows from Theorem~\ref{th:KPS}. The last maximal operator is over the complement of the
small arc containing the singularity of the weight. Hence there is no influence of the weight
on this operator and its boundedness follows trivially from Theorem~\ref{th:KPS}.
%%%%%%%%%%%%%%%%%%%%%%%%%%%%%%%%%%%%%%%%%%%%%%%%%%%%%%%%%%%%%%%%%%%%%%%%%%%%%%%%%%%%%%%%
\section{Sufficient condition involving Muckenhoupt weights}
\label{sec:Muckenhoupt-ersatz}
Although a complete characterization of weights for which $M$ is bounded
on weighted variable Lebesgue spaces is still unknown, one of the most significant
recent results to achieve this aim is the following sufficient condition
(see \cite[Theorem~${\rm A}^\prime$]{KSS07-JFSA}).
%%%%%%%%%%%%%%%%%%%%%%%%%%%%%%%%%%%%%%%%%%%%%%%%%%%%%%%%%%%%%%%%%%%%%%%%%%%%%%%%%%%%%%%%
\begin{thm}[V.~Kokilashvili, N.~Samko, S.~Samko]
\label{th:KSS}
Let $\Gamma$ be a Carleson curve, $p:\Gamma\to(1,\infty)$ be a continuous function
satisfying the Dini-Lipschitz condition, and $w:\Gamma\to[0,\infty]$ be a weight
such that $w^{p/p_*}\in A_{p_*}(\Gamma)$, where
%%%
\begin{equation}\label{eq:p-min}
p_*:=p_*(\Gamma):=\min_{\tau\in\Gamma}p(\tau).
\end{equation}
%%%
Then $M$ is bounded on $L^{p(\cdot)}(\Gamma,w)$.
\end{thm}
%%%%%%%%%%%%%%%%%%%%%%%%%%%%%%%%%%%%%%%%%%%%%%%%%%%%%%%%%%%%%%%%%%%%%%%%%%%%%%%%%%%%%%%%
This theorem does not contain the sufficiency portion of Theorem~\ref{th:KPS} whenever
$p$ is variable because for the weight $\varrho(\tau)=|\tau-t|^\lambda$ the condition
$\varrho^{p/p_*}\in A_{p_*}(\Gamma)$ is equivalent to $-1/p(t)<\lambda<(p_*-1)/p(t)$,
while the ``correct" interval for $\lambda$ is wider: $-1/p(t)<\lambda<(p(t)-1)/p(t)$.
This means that conditions of Theorem~\ref{th:KSS} cannot be necessary unless
$p$ is constant.

Now we apply Theorem~\ref{th:KSS} to the weight $\varphi_{t,\gamma}$.
%%%%%%%%%%%%%%%%%%%%%%%%%%%%%%%%%%%%%%%%%%%%%%%%%%%%%%%%%%%%%%%%%%%%%%%%%%%%%%%%%%%%%%%%
\begin{lem}\label{le:ersatz}
Let $\Gamma$ be a Carleson curve and $p:\Gamma\to(1,\infty)$ be a continuous function
satisfying the Dini-Lipschitz condition. If $t\in\Gamma$, $\gamma\in\C$, and
\begin{eqnarray}
0&<&
\frac{1}{p(t)}+{\rm Re}\,\gamma+
\min\{\delta_t^-{\rm Im}\,\gamma,\delta_t^+{\rm Im}\,\gamma\}
\nonumber
\\
&\le&
\frac{1}{p(t)}+{\rm Re}\,\gamma+
\max\{\delta_t^-{\rm Im}\,\gamma,\delta_t^+{\rm Im}\,\gamma\}<\frac{p_*}{p(t)},
\label{eq:ersatz-1}
\end{eqnarray}
where $p_*$ is defined by {\rm(\ref{eq:p-min})},
then $M$ is bounded on $L^{p(\cdot)}(\Gamma,\varphi_{t,\gamma})$.
\end{lem}
%%%%%%%%%%%%%%%%%%%%%%%%%%%%%%%%%%%%%%%%%%%%%%%%%%%%%%%%%%%%%%%%%%%%%%%%%%%%%%%%%%%%%%%%
\begin{pf}
Inequalities (\ref{eq:ersatz-1}) are equivalent to
%%%
\begin{eqnarray*}
0&<&
\frac{1}{p_*}+{\rm Re}\,\gamma\frac{p(t)}{p_*}+
\min\left\{
\delta_t^-{\rm Im}\,\gamma\frac{p(t)}{p_*},\delta_t^+{\rm Im}\,\gamma\frac{p(t)}{p_*}
\right\}
\nonumber
\\
&\le&
\frac{1}{p_*}+{\rm Re}\,\gamma\frac{p(t)}{p_*}+
\max\left\{
\delta_t^-{\rm Im}\,\gamma\frac{p(t)}{p_*},\delta_t^+{\rm Im}\,\gamma\frac{p(t)}{p_*}
\right\}<1.
\end{eqnarray*}
%%%
By Theorem~\ref{th:BK}, the latter inequalities are equivalent to
$\varphi_{t,\gamma p(t)/p_*}\in A_{p_*}(\Gamma)$.

Observe that that the weights $\varphi_{t,\gamma p(t)/p_*}$ and
$(\varphi_{t,\gamma})^{p/p_*}$ are equivalent and therefore belong to $A_{p_*}(\Gamma)$
only simultaneously. Indeed, from Proposition~\ref{pr:Carleson} (c) and
(\ref{eq:Dini-Lipschitz}) it follows that
%%%
\begin{eqnarray*}
&&
\frac{[\varphi_{t,\gamma}(\tau)]^{p(\tau)/p_*}}{\varphi_{t,\gamma p(t)/p_*}(\tau)}
=
\frac{\displaystyle\exp\left\{
\big({\rm Re}\,\gamma\log|\tau-t|-{\rm Im}\,\gamma\arg(\tau-t)\big)\frac{p(\tau)}{p_*}
\right\}}
{\displaystyle\exp\left\{
\frac{{\rm Re}\,\gamma\,p(t)}{p_*}\log|\tau-t|-\frac{{\rm Im}\,\gamma\, p(t)}{p_*}\arg(\tau-t)
\right\}}
\\
&&=
\exp\left\{
\left(\frac{{\rm Re}\,\gamma}{p_*}\log|\tau-t|-\frac{{\rm Im}\,\gamma}{p_*}\arg(\tau-t)\right)
\big(p(\tau)-p(t)\big)
\right\}
\\
&&=
\exp\left\{
\left(\frac{{\rm Re}\,\gamma}{p_*}\log|\tau-t|+\frac{{\rm Im}\,\gamma}{p_*}O(\log|\tau-t|)\right)
O\left(\frac{1}{-\log|\tau-t|}\right)
\right\}
\\
&&=\exp\{O(1)\}
\end{eqnarray*}
as $\tau\to t$. This immediately implies that the weights $(\varphi_{t,\gamma})^{p/p_*}$
and $\varphi_{t,\gamma p(t)/p_*}$ are equivalent because they are continuous on
$\Gamma\setminus\{t\}$.

Finally, applying Theorem~\ref{th:KSS}, we obtain that $M$ is bounded
on $L^{p(\cdot)}(\Gamma,\varphi_{t,\gamma})$. \qed
\end{pf}
%%%%%%%%%%%%%%%%%%%%%%%%%%%%%%%%%%%%%%%%%%%%%%%%%%%%%%%%%%%%%%%%%%%%%%%%%%%%%%%%%%%%%%%%
\section{Estimates of weights with one singularity by power weights}
\label{sec:estimate}
Recall that there are more convenient formulas for calculation of indices
of a regular submultiplicative function.
%%%%%%%%%%%%%%%%%%%%%%%%%%%%%%%%%%%%%%%%%%%%%%%%%%%%%%%%%%%%%%%%%%%%%%%%%%%%%%%%%%%%%%%%
\begin{thm}\label{th:submult-properties}
If $\varrho:(0,\infty)\to(0,\infty)$ is regular and submultiplicative, then
\begin{enumerate}
\item[{\rm (a)}]
\[
\alpha(\varrho)=\lim_{x\to 0}\frac{\log\varrho(x)}{\log x},
\quad
\beta(\varrho)=\lim_{x\to\infty}\frac{\log\varrho(x)}{\log x},
\]
and
\[
-\infty<\alpha(\varrho)\le\beta(\varrho)<+\infty;
\]

\item[{\rm (b)}]
$\varrho(x)\ge x^{\alpha(\varrho)}$ for all $x\in(0,1)$
and $\varrho(x)\ge x^{\beta(\varrho)}$ for all $x\in(1,\infty)$;

\item[{\rm (c)}]
given any $\eps>0$, there exists an $x_0>1$ such that
$\varrho(x)\le x^{\alpha(\varrho)-\eps}$ for all $x\in(0,x_0^{-1})$
and $\varrho(x)\le x^{\beta(\varrho)+\eps}$ for all $x\in(x_0,\infty)$.
\end{enumerate}
\end{thm}
%%%%%%%%%%%%%%%%%%%%%%%%%%%%%%%%%%%%%%%%%%%%%%%%%%%%%%%%%%%%%%%%%%%%%%%%%%%%%%%%%%%%%%%%
Part (a) is proved, for instance, in \cite[Theorem~1.13]{BK97}. Parts (b) and (c) follow
from part (a), see e.g. \cite[Corollary~1.14]{BK97}.

Fix $t_0\in\Gamma$. Let $\omega(t_0,\delta)$ denote the open arc on $\Gamma$
which contains $t_0$ and whose endpoints lie on the circle
$\{\tau\in\C:|\tau-t_0|=\delta\}$. It is clear that
$\omega(t_0,\delta)\subset\Gamma(t_0,\delta)$, however, it may happen that
$\omega(t_0,\delta)\ne\Gamma(t_0,\delta)$.
%%%%%%%%%%%%%%%%%%%%%%%%%%%%%%%%%%%%%%%%%%%%%%%%%%%%%%%%%%%%%%%%%%%%%%%%%%%%%%%%%%%%%%%%
\begin{lem}\label{le:estimate}
Let $\Gamma$ be a Carleson curve and $t_0\in\Gamma$. Suppose
$w:\Gamma\setminus\{t_0\}\to(0,\infty)$ is a continuous function and $W_{t_0}w$
is regular. Let $\eps>0$ and $\delta$ be such that $0<\delta<d_{t_0}$. Then
there exist positive constants $C_j=C_j(\eps,\delta,w)$, where $j=1,2$, such that
%%%
\begin{equation}\label{eq:estimate-1}
\frac{w(t)}{w(\tau)}
\le
C_1\left|\frac{t-t_0}{\tau-t_0}\right|^{\beta(W_{t_0}w)+\eps}
\end{equation}
%%%
for all $t\in\Gamma\setminus\omega(t_0,\delta)$ and all $\tau\in\omega(t_0,\delta)$;
and
%%%
\begin{equation}\label{eq:estimate-2}
\frac{w(t)}{w(\tau)}
\le
C_2\left|\frac{t-t_0}{\tau-t_0}\right|^{\alpha(W_{t_0}w)-\eps}
\end{equation}
%%%
for all $t\in\omega(t_0,\delta)$ and all $\tau\in\Gamma\setminus\omega(t_0,\delta)$.
\end{lem}
%%%%%%%%%%%%%%%%%%%%%%%%%%%%%%%%%%%%%%%%%%%%%%%%%%%%%%%%%%%%%%%%%%%%%%%%%%%%%%%%%%%%%%%%
\begin{pf}
Let us denote $\beta:=\beta(W_{t_0}w)$. By Theorem~\ref{th:submult-properties}(c),
for every $\eps>0$ there exists an $x_0\in(1,\infty)$ such that
%%%
\[
(W_{t_0}w)(x)\le x^{\beta+\eps}
\quad\mbox{for all}\quad x\in(x_0,\infty).
\]
From this inequality and the definition of $W_{t_0}w$ it follows that if $0<R\le d_{t_0}$
and $x\in(x_0,\infty)$, then
\[
\max_{|t-t_0|=R}w(t)
\le
x^{\beta+\eps}\min_{|\tau-t_0|=x^{-1}R}w(\tau)
=
\left(\frac{R}{|\tau-t_0|}\right)^{\beta+\eps}\min_{|\tau-t_0|=x^{-1}R}w(\tau).
\]
Hence
%%%
\begin{equation}\label{eq:estimate-3}
w(t)\le\left|\frac{t-t_0}{\tau-t_0}\right|^{\beta+\eps}w(\tau)
\end{equation}
%%%
for all $t\in\Gamma\setminus\{t_0\}$ and all $\tau\in\Gamma$ such that
$|t-t_0|/|\tau-t_0|\in(x_0,\infty)$. Put
\[
\Delta_{t_0}:=\min_{t\in\Gamma\setminus\omega(t_0,\delta)}|t-t_0|.
\]
It is clear that if $\tau\in\omega(t_0,\Delta_{t_0}/x_0)$ and $t\in\Gamma\setminus\omega(t_0,\delta)$,
then (\ref{eq:estimate-3}) holds.

Since the function
\[
f(\tau):=\frac{w(\tau)}{|\tau-t_0|^{\beta+\eps}}
\]
is continuous on $\Gamma\setminus\{t_0\}$, we have
\[
0<M_1:=\inf_{\tau\in\overline{\omega(t_0,\delta)\setminus\omega(t_0,\Delta_{t_0}/x_0)}}f(\tau),
\quad
M_2:=\sup_{\tau\in\overline{\Gamma\setminus\omega(t_0,\delta)}}f(\tau)<\infty.
\]
Hence
\[
w(t)\le M_2|t-t_0|^{\beta+\eps}
\]
for all $t\in\Gamma\setminus\omega(t_0,\delta)$ and
\[
\frac{1}{w(\tau)}\le\frac{1}{M_1|\tau-t_0|^{\beta+\eps}}
\]
for all $\tau\in\omega(t_0,\delta)\setminus\omega(t_0,\Delta_{t_0}/x_0)$.
Multiplying these inequalities, we obtain
\[
\frac{w(t)}{w(\tau)}\le\frac{M_2}{M_1}\left|\frac{t-t_0}{\tau-t_0}\right|^{\beta+\eps}
\]
for all $t\in\Gamma\setminus\omega(t_0,\delta)$ and all
$\tau\in\omega(t_0,\delta)\setminus\omega(t_0,\Delta_{t_0}/x_0)$.

Thus (\ref{eq:estimate-1}) holds for $t\in\Gamma\setminus\omega(t_0,\delta)$ and all
\[
\tau\in\omega(t_0,\Delta_{t_0}/x_0)\cup
[\omega(t_0,\delta)\setminus\omega(t_0,\Delta_{t_0}/x_0)]=\omega(t_0,\delta)
\]
with $C_1:=\max\{1,M_2/M_1\}$. Estimate (\ref{eq:estimate-2}) is proved by analogy.
\qed
\end{pf}
%%%%%%%%%%%%%%%%%%%%%%%%%%%%%%%%%%%%%%%%%%%%%%%%%%%%%%%%%%%%%%%%%%%%%%%%%%%%%%%%%%%%%%%%
\section{Proof of Theorem~\ref{th:main}}
\label{sec:proof}
The idea of the proof is borrowed from \cite[Theorem~B]{KSS07-JFSA}.
Fix $t_0\in\Gamma$ and $\gamma\in\C$. Notice that we omitted the subscript for $t_0$ in the
formulation of the theorem for brevity. It is easily seen that $M$ is bounded on
$L^{p(\cdot)}(\Gamma,\varphi_{t_0,\gamma})$ if and only if the operator
\[
(M_{t_0,\gamma}f)(t):=\sup_{\eps>0}\frac{\varphi_{t_0,\gamma}(t)}{|\Gamma(t,\eps)|}
\int_{\Gamma(t,\eps)}\frac{|f(\tau)|}{\varphi_{t_0,\gamma}(\tau)}\,|d\tau|
\quad
(t\in\Gamma)
\]
is bounded on $L^{p(\cdot)}(\Gamma)$.

As it was already mentioned in the introduction, the function $W_{t_0}\varphi_{t_0,\gamma}$
is regular and submultiplicative for every $\gamma\in\C$ and
%%%
\begin{eqnarray*}
\alpha
&:=&
\alpha(W_{t_0}\varphi_{t_0,\gamma})
=
{\rm Re}\,\gamma+
\min\{\delta_{t_0}^-{\rm Im}\,\gamma,\delta_{t_0}^+{\rm Im}\,\gamma\},
\\
\beta
&:=&
\beta(W_{t_0}\varphi_{t_0,\gamma})
=
{\rm Re}\,\gamma+
\max\{\delta_{t_0}^-{\rm Im}\,\gamma,\delta_{t_0}^+{\rm Im}\,\gamma\}.
\end{eqnarray*}
%%%
With these notations, conditions (\ref{eq:main-conditions}) have the form
\[
0<\frac{1}{p(t_0)}+\alpha,
\quad
\frac{1}{p(t_0)}+\beta<1.
\]
In this case there is a small $\eps>0$ such that
%%%
\begin{equation}\label{eq:main-1}
0<\frac{1}{p(t_0)}+\alpha-\eps
\le
\frac{1}{p(t_0)}+\beta+\eps<1.
\end{equation}
%%%
Since $p:\Gamma\to(1,\infty)$ is continuous and $1/p(t_0)+\beta<1$, we can choose a number $\delta\in(0,d_{t_0})$ such
that the arc $\omega(t_0,\delta)$, which contains $t_0$ and has the endpoints on the circle
$\{\tau\in\C:|\tau-t_0|=\delta\}$, is so small that $1+\beta p(t_0)<p_*$, where
\[
p_*:=p_*(\omega(t_0,\delta))=\min_{\tau\in\overline{\omega(t_0,\delta)}}p(\tau).
\]
Hence
%%%
\begin{equation}\label{eq:main-2}
0<\frac{1}{p(t_0)}+\alpha\le\frac{1}{p(t_0)}+\beta<\frac{p_*}{p(t_0)}.
\end{equation}
%%%
Let us denote by $\chi_\Omega$ the characteristic function of a set $\Omega\subset\Gamma$.
For $f\in L^{p(\cdot)}(\Gamma)$, we have
%%%
\begin{eqnarray}
M_{t_0,\gamma}f
&\le&
\chi_{\omega(t_0,\delta)}M_{t_0,\gamma}\chi_{\omega(t_0,\delta)}f+
\chi_{\Gamma\setminus\omega(t_0,\delta)}M_{t_0,\gamma}\chi_{\omega(t_0,\delta)}f
\nonumber
\\
&&+
\chi_{\omega(t_0,\delta)}M_{t_0,\gamma}\chi_{\Gamma\setminus\omega(t_0,\delta)}f+
\chi_{\Gamma\setminus\omega(t_0,\delta)}M_{t_0,\gamma}\chi_{\Gamma\setminus\omega(t_0,\delta)}f.
\label{eq:main-3}
\end{eqnarray}
%%%

From (\ref{eq:main-2}) and Lemma~\ref{le:ersatz} we conclude that $M_{t_0,\gamma}$ is bounded on
$L^{p(\cdot)}(\omega(t_0,\delta))$. Hence $\chi_{\omega(t_0,\delta)}M_{t_0,\gamma}\chi_{\omega(t_0,\delta)}I$
is bounded on $L^{p(\cdot)}(\Gamma)$.

From Lemma~\ref{le:estimate} it follows that
%%%
\begin{equation}\label{eq:main-4}
\chi_{\Gamma\setminus\omega(t_0,\delta)}M_{t_0,\gamma}\chi_{\omega(t_0,\delta)}f
\le
C_1
\chi_{\Gamma\setminus\omega(t_0,\delta)}M_{t_0,\beta+\eps}\chi_{\omega(t_0,\delta)}f
\le
C_1M_{t_0,\beta+\eps}f
\end{equation}
%%%
and
%%%
\begin{equation}\label{eq:main-5}
\chi_{\omega(t_0,\delta)}M_{t_0,\gamma}\chi_{\Gamma\setminus\omega(t_0,\delta)}f
\le
C_2
\chi_{\omega(t_0,\delta)}M_{t_0,\alpha-\eps}\chi_{\Gamma\setminus\omega(t_0,\delta)}f
\le
C_2M_{t_0,\alpha-\eps}f,
\end{equation}
%%%
where $C_1$ and $C_2$ are positive constants depending only on $\eps,\delta,\gamma$, and $t_0$.
Inequalities (\ref{eq:main-1}), Theorem~\ref{th:KPS}, and inequalities
(\ref{eq:main-4})--(\ref{eq:main-5}) imply that the operators
$\chi_{\Gamma\setminus\omega(t_0,\delta)}M_{t_0,\gamma}\chi_{\omega(t_0,\delta)}I$
and
$\chi_{\omega(t_0,\delta)}M_{t_0,\gamma}\chi_{\Gamma\setminus\omega(t_0,\delta)}I$
are bounded on $L^{p(\cdot)}(\Gamma)$.

Finally, since $\Gamma\setminus\omega(t_0,\delta)$ does not contain the singularity
of the weight $\varphi_{t_0,\gamma}$ which is continuous on $\Gamma\setminus\{t_0\}$, there
exists a constant $C_3>0$ such that
\[
\chi_{\Gamma\setminus\omega(t_0,\delta)}M_{t_0,\gamma}\chi_{\Gamma\setminus\omega(t_0,\delta)}f
\le
C_3Mf.
\]
Then Theorem~\ref{th:KPS} and the above estimate yield the boundedness of the operator
$\chi_{\Gamma\setminus\omega(t_0,\delta)}M_{t_0,\gamma}\chi_{\Gamma\setminus\omega(t_0,\delta)}I$
on $L^{p(\cdot)}(\Gamma)$.

Thus all operators on the right-hand side of (\ref{eq:main-3}) are bounded on $L^{p(\cdot)}(\Gamma)$.
Therefore the operator on the left-hand side of (\ref{eq:main-3}) is bounded, too. This
completes the proof.
\qed
%%%%%%%%%%%%%%%%%%%%%%%%%%%%%%%%%%%%%%%%%%%%%%%%%%%%%%%%%%%%%%%%%%%%%%%%%%%%%%%%%%%%%%%%

\end{document}